\documentclass[preprint,12pt]{elsarticle}
\journal{Applied Numerical Mathematics}

\usepackage{multicol,amsmath,amssymb,latexsym,amsfonts}
\usepackage{paralist}

\newcommand{\x}{{\bf x}}
\newcommand{\n}{{\bf n}}
\newcommand{\y}{{\bf y}}

\newcommand{\pderiv}[2]{\frac{\partial #1}{\partial #2}}

\newcommand{\eqr}[1]{~(\ref{#1})}
\newcommand{\figr}[1]{figure~\ref{#1}}

\newcommand{\Tabr}[1]{Table~\ref{#1}}

 \usepackage{graphicx}
 \usepackage{subfigure}

\journal{Applied Numerical Mathematics}
\begin{document}




\begin{frontmatter}

\title{Fast integral equation methods for the heat equation and the modified Helmholtz equation in two dimensions}
 \author{M.C.A. Kropinski\corref{nserc}\fnref{mcak}}
 \author{Bryan Quaife\fnref{mcak}}
 \address[mcak]{Department of Mathematics, Simon Fraser University,
 Burnaby, British Columbia, Canada V5A 1S6}
 \cortext[nserc]{Supported in part by Natural Sciences and Engineering Research Council of Canada Grant RGPIN 203326}

\begin{abstract}
We present an efficient integral equation approach to solve the heat equation, $u_t (\x) - \Delta
  u(\x) = F(\x,t)$, in a two-dimensional, multiply connected domain, and with Dirichlet boundary conditions. 
Instead of using integral equations based on the heat kernel, we take the approach of discretizing in time, first. 
This leads to a non-homogeneous modified Helmholtz equation that is solved at each time step. The solution to this equation is formulated as a volume potential plus a double layer potential.
The volume potential is evaluated using a fast multipole-accelerated solver.
The boundary conditions are then satisfied by solving an integral equation for the homogeneous modified Helmholtz equation.  
The integral equation solver is also accelerated by the fast multipole method (FMM). 
For a total of $N$ points in the discretization of the boundary and the domain, the total computational cost per time step is $O(N)$ or $O(N\log N)$. 
\end{abstract}

\end{frontmatter}

\section{Introduction}
Integral equation methods offer an attractive alternative to conventional finite difference and finite elements for solving partial differential equations that arise in science and engineering. 
They offer several advantages: complex physical boundaries are easy to incorporate, the ill-conditioning associated with directly discretizing the governing equation is avoided, high-order accuracy is easier to attain, and far-field boundary conditions are handled naturally.
Nevertheless, integral equation methods have been slow to be adopted for large-scale problems, and the  reason for this is clear. 
For non-homogeneous problems, unknowns are distributed at all points of the domain. 
On an $N\times N$ mesh, integral equation formulations lead to dense $N^2 \times N^2$ matrices.
This situation is further exacerbated if the problem is time dependent. 

In recent years, fast algorithms have been developed for integral equations in a variety of settings. 
An early example is the work by Greenbaum et al. in \cite{mult:conn} which formulates integral equation methods for Laplace's equation in multiply connected domains.
The solution to these integral equations was accelerated by the original version Greengard and Rokhlin's fast multipole method (FMM) \cite{fmm,gr}.
For $N$ points in the discretization of the domain, the solution to the integral equation requires only $O(N)$ operations. 
Subsequently, similar tools have been developed for other elliptic partial differential equations, such as Poisson's equation \cite{greengard:ethridge}, the Stokes equations \cite{stokes:flow}, and the modified Helmholtz equation\cite{modified:helmholtz,KROP_QUAIFE_1}. 
Nishimura provides a comprehensive review of fast multipole accelerated integral equation methods in \cite{nishimura}. 

For time-dependent problems, there is added complexity. 
There is recent work by Li and Greengard in \cite{heat_solver_1,heat_solver_2} on developing fast solvers for the heat equation. 
Their methods use an integral equation formulation based on the heat kernel, and while robust, are highly complicated. 
To date, the full set of tools required for solving the heat equation in interior or exterior domains with complex boundaries has not been constructed. 

We present an alternative approach to solving the heat equation and other time-dependent partial differential equations that rely on available fast algorithms for elliptic problems. 
The temporal discretization of the heat equation gives rise to  the modified Helmholtz, or linearized Poisson-Boltzmann, equation
\begin{equation}
  u(\x) - \alpha^2 \Delta u(\x) = B(\x,t),
  \label{eq:mod_lap_hom}
\end{equation}
where $\alpha^2 = O(\Delta t)$. 
This equation seems to have received little attention in the literature, particularly in this context. 
Hsiao and Maccamy formulate integral equations of the first kind in \cite{hsiao} for this equation. 
In \cite{modified:helmholtz}, Cheng et al. present a
fast direct solver for\eqr{eq:mod_lap_hom} in two dimensions on the unit square.  
The solution is expressed as a volume potential, and the direct solver is accelerated using a new version of the fast multipole method \cite{screened_coulomb, new_FMM}. 
The solver is fully adaptive and the computational costs are comparable to those of FFT based methods.
In \cite{KROP_QUAIFE_1}, we present fast, well-conditioned integral equation methods to solve the homogeneous equation in unbounded and bounded multiply-connect domains, with Dirichlet or Neumann boundary conditions. 
What we present here is a preliminary study on coupling the methods discussed in \cite{modified:helmholtz} to those discussed in \cite{KROP_QUAIFE_1} in order to solve practical problems in complicated domains. 

The paper is organized as follows: In Section 2, we discuss the temporal discretization of the heat equation and how this gives rise to the modified Helmholtz equation. 
We outline the corresponding potential theory; decouple the solution into a volume potential and double layer potential with an unknown density, and then formulate an integral equation for the Dirichlet problem in a bounded multiply connected domain. In Section 3,
numerical methods for the evaluation of the volume potential and the integral equation are 
discussed. 
Numerical examples are presented in Section 5.

\begin{figure}[t]
     \centering
$\begin{array}{c}
\includegraphics[height=2in]{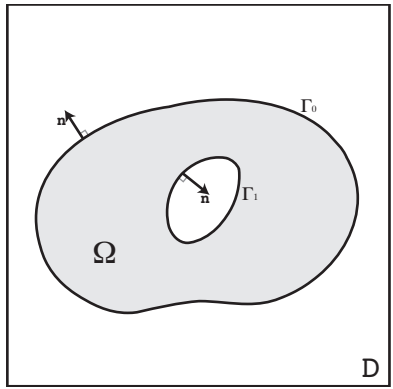}
\end{array}$
  \caption{\em The physical domain $\Omega$ (in gray) is embedded in the unit square $D$.  The physical outer boundary is denoted by $\Gamma_0$ and the interior boundary is $\Gamma_1$.}
  \label{fig1}
\end{figure}
\section{Potential theory and the integral equation formulation}
Consider the isotropic inhomogeneous heat equation in a bounded domain $\Omega$ in $\mathbb{R}^2$:
\begin{equation}
  u_t(\x,t) -\Delta u(\x,t) = F(\x,t),  \qquad  t>0,  \label{pde} 
\end{equation}
with Dirichlet boundary conditions:
\[
  u(\x,t) =  f(\x,t), \qquad  \x \in \Gamma, \; t>0, 
\]
and initial conditions
\[
  u(\x,0) = u_0(\x), \qquad  {\x} \in \Omega.  
\]
One approach to solving this boundary value problem would be base an integral equation approach on the fundamental solution for the heat equation. 
Fast algorithms are developed for just such an approach in \cite{heat_solver_1,heat_solver_2}. 
We approach from  a different starting point by discretizing\eqr{pde} in time, first, and then reformulate in terms of an integral equation.

In order to prevent a severe time-step restriction, linearly implicit or IMEX schemes \cite{Ruuth} are generally used for marching in time. 
In such schemes, the diffusive term is treated implicitly, while the remainder is treated explicitly.
Regardless of the details of the particular choice of IMEX scheme, the temporal discretization of\eqr{pde} yields
\begin{equation}
\begin{array}{rcll}
  u^{N+1} - \alpha^2 \Delta u^{N+1} & = &  B(u^N,u^{N-1},\cdots), & \qquad {\x} \in \Omega,  \smallskip \\
  u^{N+1} & =  & f(\x,t), &  \qquad \x \in \Gamma, \smallskip \\
  u^0 & =  & u_0(\x), & \qquad \mbox{at } t = 0,
 \end{array}
\label{mod_helm_forced} 
\end{equation}
where $t=(N+1) \Delta t$ is the current time.
The simplest such scheme is the first-order backward Euler method, which is of the form\eqr{mod_helm_forced} with   
\[
   \alpha^2 =  \Delta t, \qquad B = u^N + \Delta t F^N. 
\]
A second-order method is the extrapolated Gear method, for which
\[
\alpha^2 =  \frac{2}{3} \Delta t, \qquad B = \frac{4}{3} u^N -\frac{1}{3} u^{N-1}+
                 \frac{4}{3}  \Delta t F^N  - \frac{2}{3} \Delta t F^{N-1}. \label{exg}
\]

In order to solve\eqr{mod_helm_forced} by means of integral equations, we first represent the solution at each time step as
\[
U(\x) = U^p(\x) + U^h(\x),  
\]
where
\begin{equation}
   U^p - \alpha^2 \Delta U^p = B, \qquad    \x \in \Omega, 
   \label{particular}
\end{equation}
and $U^h$ satisfies
\begin{equation}
   \begin{array}{rcll}
   U^h  - \alpha^2 \Delta U^h & = & 0, &   \qquad \x \in \Omega, \smallskip \\
   U^h & =  & g(\x) , &   \qquad \x \in \Gamma,  
   \end{array}
 \label{mod_helm_hom}
\end{equation}
where $g(\x) =  f(\x) - U^p(\x)$, $\x \in \Gamma$.

The solution to\eqr{particular} is represented as 
\begin{equation}
   U^p(\x) = \int_\Omega B(\y) G(\y - \x) \, dA_{\y}, \label{vol_potential}
\end{equation}
where 
the free space Green's function $G(\x) $ for the operator $1-\alpha^2\Delta$ is given by
the zeroth-order modified Bessel function of the second
kind,
\[
  G(\x) = \frac{1}{2\pi\alpha^2} K_0 \left( \frac{|\x|}{\alpha} \right) .
\]
We seek the solution to\eqr{mod_helm_hom} in the form of a double layer potential
\[
  U^h(\x)=\frac{1}{2\pi\alpha^2}\int_{\Gamma}\pderiv{}{\nu_{\y}}K_{0}
  \left(\frac{|\y-\x|}{\alpha} \right) \sigma(\y) ds(\y), \qquad \x\in\Omega
\]
where  $\sigma(\y)$ is the value of an unknown
density at the boundary point $\y$, and $\partial / \partial \nu_\y$
represents the outward normal derivative at the point $\y$.
The unknown density $\sigma$ is found by solving an integral equation derived to ensure the boundary conditions in\eqr{mod_helm_hom} are satisfied. 

In \cite{KROP_QUAIFE_1}, we derive integral equation formulations to solve\eqr{mod_helm_hom} in bounded and unbounded multiply-connected domains. 
Here, we summarize the results for the bounded case.  
Observing that $K_0 (z) \sim -\log(z)$ as $z \rightarrow 0$ \cite{ABRAM}, 
we can determine the jump relations of the double layer potential as these are well known for the logarithmic kernel. 
Thus for  any point $\x$ on the boundary $\Gamma$,
\begin{eqnarray*}
  \lim_{\tiny{\begin{array}{l}
        \x'  \rightarrow  \x  \\
        \x'  \in  \Omega 
      \end{array} }}
  u(\x')& = &  - \frac{1}{2\alpha^2}\sigma(\x)+
  \frac{1}{2\pi\alpha^2}\int_{\Gamma}\pderiv{}{\nu_{\y}}
  K_{0}\left(\frac{|\y-\x|}{\alpha}\right)\sigma(\y)ds(\y).
\end{eqnarray*}
Substituting the above into the boundary condition in\eqr{mod_helm_hom} results in the integral equation 
\begin{equation}
  \label{dirichlet:inteqn}
 -\frac{1}{2\alpha^2}\sigma(\x)+\frac{1}{2\pi\alpha^2}
  \int_{\Gamma}\pderiv{}{\nu_{\y}}K_{0}\left(\frac{|\y - \x|}{\alpha}\right)
  \sigma(\y) \, ds(\y) = g(\x) .
\end{equation}
The kernel in the above is continuous along
$\Gamma$, 
\[
  \lim_{\tiny{\begin{array}{r}\y \rightarrow \x \smallskip \\
        \x, \y \in \Gamma
      \end{array} }}
  \pderiv{}{\nu_{\y}}K_{0}\left(\frac{|\y-\x|}{\alpha}\right) = 
  -\frac{1}{2}\kappa(\x),
\]
where $\kappa(\x)$ denotes the curvature of $\Gamma$ at the point
$\x$,
(here, we have again used the logarithmic behaviour of $K_0(z)$ as $z \rightarrow 0$).
In addition,\eqr{dirichlet:inteqn} has no nontrivial
homogeneous solutions. Therefore, by the Fredholm alternative,\eqr{dirichlet:inteqn} has a unique solution for any integrable data
$f(\x)$.
\section{Numerical methods}
We now discuss the numerical methods for solving the modified Helmholtz equation\eqr{mod_helm_forced}.
In Section 3.1, we first outline the methods involved in evaluating the volume potential\eqr{vol_potential}. 
In Section 3.2, we next discuss the coupling of the volume potential to the double layer potential through the boundary conditions and solving the corresponding integral equation\eqr{dirichlet:inteqn}. 
The work presented in this section is discussed in more detail in \cite{KROP_QUAIFE_1}.
While the FMM plays an integral role in our numerical methods, its implementation is standard and has been discussed at length in other works. 
For details on the FMM, we refer the reader to the work in \cite{gr,fmm} for the original FMM, \cite{new_FMM} for the new FMM algorithm based on exponential expansions to represent the far-field interactions, and \cite{modified:helmholtz} for the FMM applied to the modified Helmholtz equation in two dimensions.

\subsection{The volume potential}
In \cite{modified:helmholtz}, methods for the rapid evaluation of the volume potential for the modified Helmholtz equation on the unit square were developed.
These methods use adaptive mesh refinement, and the evaluation is accelerated using the new version of the FMM. 
The total computational cost is comparable to that of an FFT-based fast solver.
In order to employ these methods, we must first extend the right-hand side of\eqr{pde}, which is defined in $\Omega$, to be defined throughout $D$. 
One naive way to define an extension $\bar{B}(\x,t)$ of $B(\x,t)$ to $D$ is the following: 
\[
   \tilde{B}(\x,t) = \left\{ \begin{array}{rl}
                                       B(\x,t), & \qquad \x \in \Omega \\
                                       0, & \qquad \x \in D \backslash  \Omega.
                                    \end{array} \right.
\]
However, the discontinuity across $\Gamma$ would likely cause excessive grid refinement  in the vicinity of the boundary, and consequently slow down the evaluation (c.f. example 4.4 in \cite{greengard:ethridge}; up to $O(10^5)$ grid points are needed to solve $\Delta u=1$ inside a disk and $\Delta u = 0$ outside). 
For multiple-component boundaries, this would simply become too expensive. 
For the purposes of this present work, we consider problems with simple boundary
conditions, only, so that $B$ can be extended continuously from $\Omega$ to
$D$.  
(We save the consideration of more general-purpose methods to extend $B$ for future work.)
For example, if we are solving the homogeneous heat equation ($F(\x,t) = 0$, then $B(\x)$ consists of a linear combination of the solution $u(\x,t)$ at previous time steps. 
If the boundary conditions are constant, $u=C_0$ on $\Gamma_0$, $u=C_1$ on $\Gamma_1$, say, then we can define a continuous $\tilde{B}(\x,t)$:
\[
   \tilde{B}(\x,t) = \left\{ \begin{array}{rl}
                                       B(\x,t), & \qquad \x \in \Omega, \\
                                       C_0, & \qquad \mbox{exterior to } \Gamma_0, \\
                                       C_1, & \qquad \mbox{exterior to } \Gamma_1.
                                    \end{array} \right.
\]
We now use the methods presented in \cite{modified:helmholtz} to evaluate
\begin{equation}
   \tilde{U}^p(\x) = \int_D\tilde{B}(\x,t) G(\x-\y) \, dA_{\y},   \label{volume_potential}
\end{equation}
and let $ U^p(\x) = \tilde{U}^p(\x)$ for $\x\in\Omega$.

The FMM uses an adaptive quad-tree structure in order to superimpose a hierarchy of refinement on the computational domain. 
The unit square $D$ is considered to be grid level 0. Grid level $l+1$ is obtained recursively by subdividing each square (or node) $s$ at level $l$ into four equal parts; these are called the ``children'' of $s$. 
Adaptivity is achieved by allowing different levels of refinement throughout the tree.
The square $s$ is refined if the error of a fourth order polynomial interpolating $\bar{B}$ is larger than some preset tolerance.  
We denote the childless nodes in the quad-tree as $D_i$, $i=1, \cdots, M$, where $M$ is the total number of such nodes. 
We assume we are given $\bar{B}$ on a cell-centered $4\times 4$ grid for each $D_i$. 
Thus, $N_D = 16 \times M$ is the total number of grid points in $D$. 
To obtain fourth-order accuracy, these 16 points are used to construct a fourth-order polynomial to $\bar{B}$ of the form
\[
  \tilde{B}(\x) \approx \sum_{j=1}^{10} c^i_j \, p_j (\x-\x^i), \qquad \x \in D_i, 
\] 
where $\x^i$ is the centre of $D_i$ and $\{ p_{j} \}$ are the standard basis functions for polynomials up to
order three
(see \cite{modified:helmholtz,greengard:ethridge} for details on the approximating polynomial).
Therefore, the evaluation of\eqr{volume_potential} is approximated by 
\[
  \tilde{U}^p (\x) \approx \sum_{i=1}^M 
    \int_{D_i} G(\x - \y) \, \sum_{j=1}^{10} c^i_j \, p_j (\x-\x^i) \, dA_{\y} .
 \]
 Direct evaluation of this potential at all $N_D$ grid points would require $O(N^2_D)$ work. 
This expense is reduced to $O(N_D)$ work, by using the FMM.

Once we have obtained the values of $\bar{U}^p$ at all $N_D$ points, we must next evaluate the volume potential on $\Gamma$ in order to impose the correct boundary conditions in\eqr{mod_helm_hom}. 
Again, we construct the approximating polynomial to $\bar{U}^p$ on each $D_i$, 
\[
  \tilde{U}^p(\x) \approx \sum_{j=1}^{10} C^i_j \, p_j (\x-\x^i), \qquad \x \in D_i.  
\]
Every point $\y \in \Gamma$ belongs to a particular childless node $D_i$. 
Determining the value of $U^p$ at $\y$ is simply a matter of evaluating the appropriate approximating polynomial at that point.

\subsection{The integral equation}
We now discuss the numerical methods to solve\eqr{dirichlet:inteqn}.
We assume each component curve $\Gamma_k$, $k=0,1$ is
parametrized by $\y^k(\alpha)$, where $\alpha \in [0, 2\pi)$. Similarly, $\sigma^k(\alpha)$ refers to the restriction of the density $\sigma$ on $\Gamma_k$.
We are given $N$ points equispaced with respect
to $\alpha$. Thus the mesh spacing is $h = 2\pi/N$, and the total
number of descretization points is $N_\Gamma =2 N$. Associated with each
such point, denoted by $\y^k_j$, is an unknown density $\sigma^k_j$.

In order to approximate the integral operator in\eqr{dirichlet:inteqn}, we use 
hybrid gauss-trapezoidal quadrature rules developed by Alpert \cite{alpert:quad:rules} which are tailored for integrands with logarithmic singularities. 
These quadratures are of order $h^p \log h$. 
The order $p$ determines the nodes $u_n$ and weights  $v_n$, $n=1,\cdots, l$, which are used for the quadrature within the interval $\alpha \in [\alpha_j- h a, \alpha_j+h a]$, on $\Gamma_k$ ($a$ is also determined by $p$). 
Outside of this interval, the quadrature is essentially the trapezoid rule.
Applying this quadrature to\eqr{dirichlet:inteqn} yields
\begin{align}
    \sigma_j^k & - \frac{h}{\pi} \left\{
    \sum_{\begin{array}{l}
                             m=0  \\
                            m \ne k
                         \end{array} }^1  \sum_{n=1}^{N} K(\y^m_n,\y^k_j)\, \sigma^m_n + \sum_{n=j+a}^{N + j-a} K(\y^k_n,\y^k_j)\, \sigma^k_n  \right\}  
                         \nonumber \\
    &   -\frac{1}{\pi} \sum_{\begin{array}{l}
                             n=-l  \\
                            n \ne 0
                         \end{array} }^l \, u_{|n|} 
             K(\y^k_{j + \frac{ n}{|n|} v_{|n|}}, \y^k_j ) 
            \, \sigma^k_{j + \frac{n}{|n|} v_{|n|}} 
             =  -2 \alpha^2  g(\y_j^k),          
             \label{discrete_sys}
\end{align}
where
\[
K(\y,\x) = \frac{1}{\alpha} K_1 \left( \frac{|\y - \x|}{\alpha} \right)
                         \frac{\y-\x }{|\y - \x|} \cdot {\bf n}_{\y}. 
\]
The outward pointing normal $\n_\y$ at each point is obtained pseudospectrally.
In the middle sum, we invoke periodicity of all functions on $\Gamma_k$, or equivalently, $j+N = j$. In the final sum, we are required to know $\sigma$ at intermediate values to the nodal values.  These are found through Fourier interpolation.

Equation\eqr{discrete_sys} is a linear system that is solved iteratively using the generalized minimum residual method GMRES \cite{SAAD}. 
The bulk of the work at each iteration lies in evaluating\eqr{discrete_sys} at the current solution update. 
If this was done directly, it would require $O(N^2_\Gamma)$ work. 
This evaluation can be reduced to $O(N_\Gamma)$, again using the FMM (c.f. \cite{KROP_QUAIFE_1} for more discussion on the details of implementation). 
Since the number of iterations needed to solve a Fredholm equation of the second kind to a fixed precision is bounded independent of $N_\Gamma$, we can estimate the total cost of solving\eqr{discrete_sys} by
\[
   I(\epsilon) C(\epsilon) N_\Gamma,
\]
where $I(\epsilon)$ is the number of GMRES iterations needed to reduce the residual error to $\epsilon$, and $C(\epsilon)$ is the constant of proportionality in the FMM. 

\section{Numerical results}
The algorithms described above have been implemented in Fortran. 
Here, we illustrate their performance on a variety of examples.
The tolerance for the residual error in GMRES is set to $10^{-10}$. 
All timings cited are for a Mac Pro 2.1 with two 3GHz Quad-Core Intel Xeon processors.

{\it{EXAMPLE 1.}}
In this example, we solve the forced heat equation on a domain for which an analytical solution is known. 
The domain $\Omega$ is bounded between $\Gamma_0$,  a circle of radius 0.4, and $\Gamma_1$, a circle of radius 0.1. Both $\Gamma_0$ and $\Gamma_1$ are centered at the origin.
The forcing term is
\[
F(\x,t) = 400 \cos(20|\x|) + 20 \frac{\sin(20|\x|)}{|\x|}
\]
The following is the exact solution to the heat equation:
\begin{align*}
  u(x,t)=e^{-\lambda^{2}t}\left[Y_{0}(0.1\lambda)J_{0}(\lambda|\x|)-
  J_{0}(0.1\lambda)Y_{0}(\lambda|\x|)\right]+\cos (20 |\x|), 
\end{align*}
where $J_0$ is the Bessel function of the first kind of order zero, and $Y_0$ is the Bessel function of the second kind of order zero. 
We have chosen $\lambda$ so that the time-dependant term vanishes on both boundaries $\lambda \approx 10.244$. 
The boundary conditions, then, are 
\[
   f(\x) = \left\{ \begin{array}{rl}
                            \cos(8), & \qquad \x \in \Gamma_0, \\
                            \cos(2), & \qquad \x \in \Gamma_1.
                         \end{array} \right.
\]
We calculate solutions to the heat equation using the IMEX Euler and SBDF methods, up to $t = 1.0e-2$. The results are summarized in \Tabr{table1}.

\begin{table}[htbp]
  \begin{center}
\begin{tabular*}{\textwidth}{@{\extracolsep{\fill}}ccc}
    \hline
  $ \Delta t$ & $\mbox{Error}_1$ & $\mbox{Error}_2$ \\  \hline
    $2.0 \times 10^{-3}$ & $1.37 \times 10^{-2}$ & $1.65 \times 10^{-3}$ \\
    $1.0 \times 10^{-3}$ &  $7.09 \times 10^{-3}$ & $4.74 \times 10^{-4}$ \\
    $5.0 \times 10^{-4}$ & $3.59 \times 10^{-3}$ & $1.23 \times 10^{-4}$ \\
    $2.5 \times 10^{-4}$ & $1.78 \times 10^{-4}$ & $3.25 \times 10^{-5}$ \\
    \hline
  \end{tabular*}
  \end{center}
\caption{Temporal Error using IMEX Euler and extrapolated Gear Method.  $\mbox{Error}_1$ indicates the error for the IMEX Euler method,  $\mbox{Error}_2$ is for extrapolated Gear. }
\label{table1}
\end{table}

\begin{figure}[htps]
     \centering
$\begin{array}{ccc}
\includegraphics[height=1.75in]{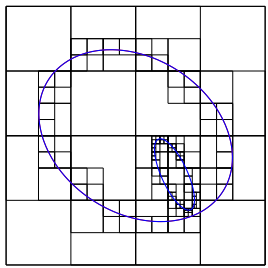} &
\end{array}$
  \caption{\em The quadtree structure for Example 2. }
  \label{figure2}
\end{figure}
\begin{figure}[htps]
     \centering
$\begin{array}{ccc}
\includegraphics[height=1.in]{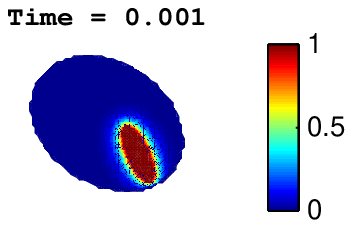} &
\includegraphics[height=1.in]{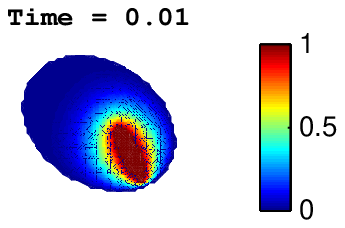} &
\includegraphics[height=1.in]{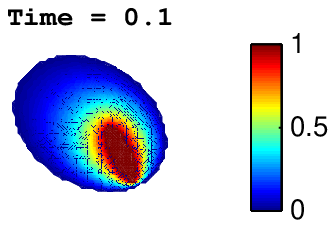} 
\end{array}$
  \caption{\em The solution to Example 2. $N_\Gamma=1024$, $N_D = 65536$. Total CPU time for 100 time steps about 30 minutes. }
  \label{figure3}
\end{figure}

{\it{EXAMPLE 2.}} In this example, we solve the homogeneous heat equation in an elliptical region containing an off-center elliptical hole. The boundary conditions are
\[
   f(\x) = \left\{ \begin{array}{rl}
                            0, & \qquad \x \in \Gamma_0, \\
                            1, & \qquad \x \in \Gamma_1.
                         \end{array} \right.
\]
The quadtree structure for the solution homogeneous equation is shown in \figr{figure2}
The results are shown in \figr{figure3}.

\begin{figure}[htps]
     \centering
$\begin{array}{cc}
\includegraphics[height=1.5in]{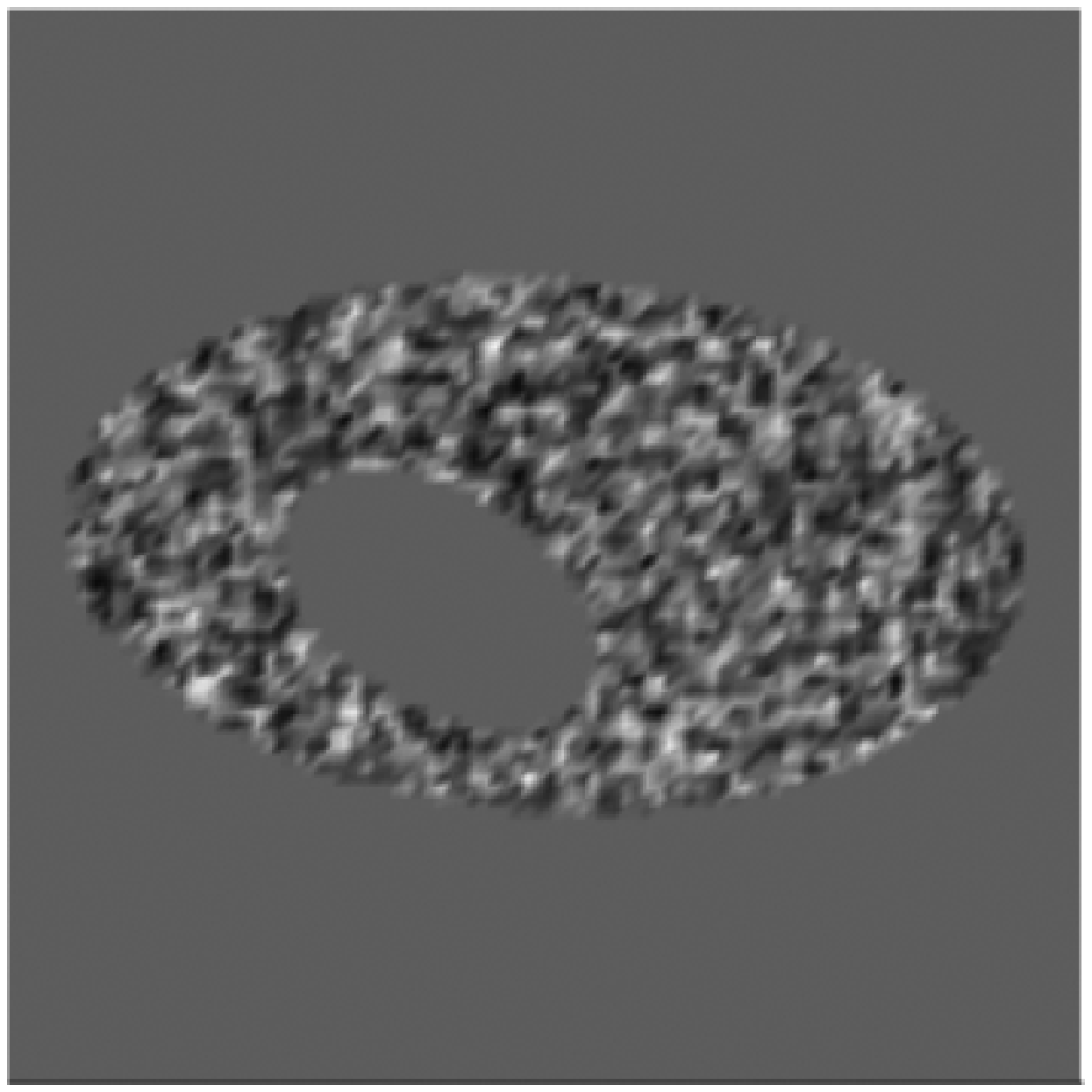} &
\includegraphics[height=1.5in]{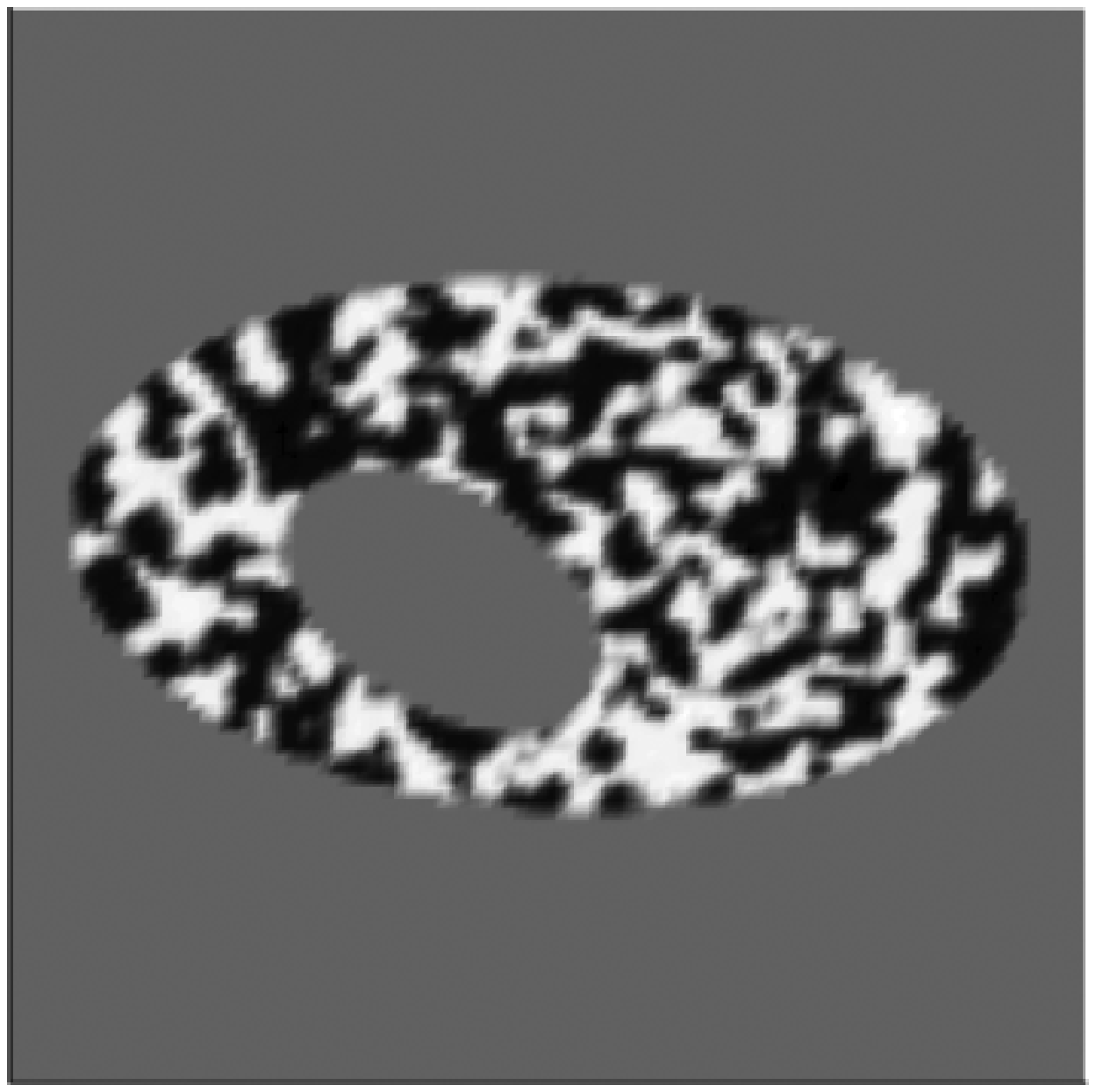} \\
\includegraphics[height=1.5in]{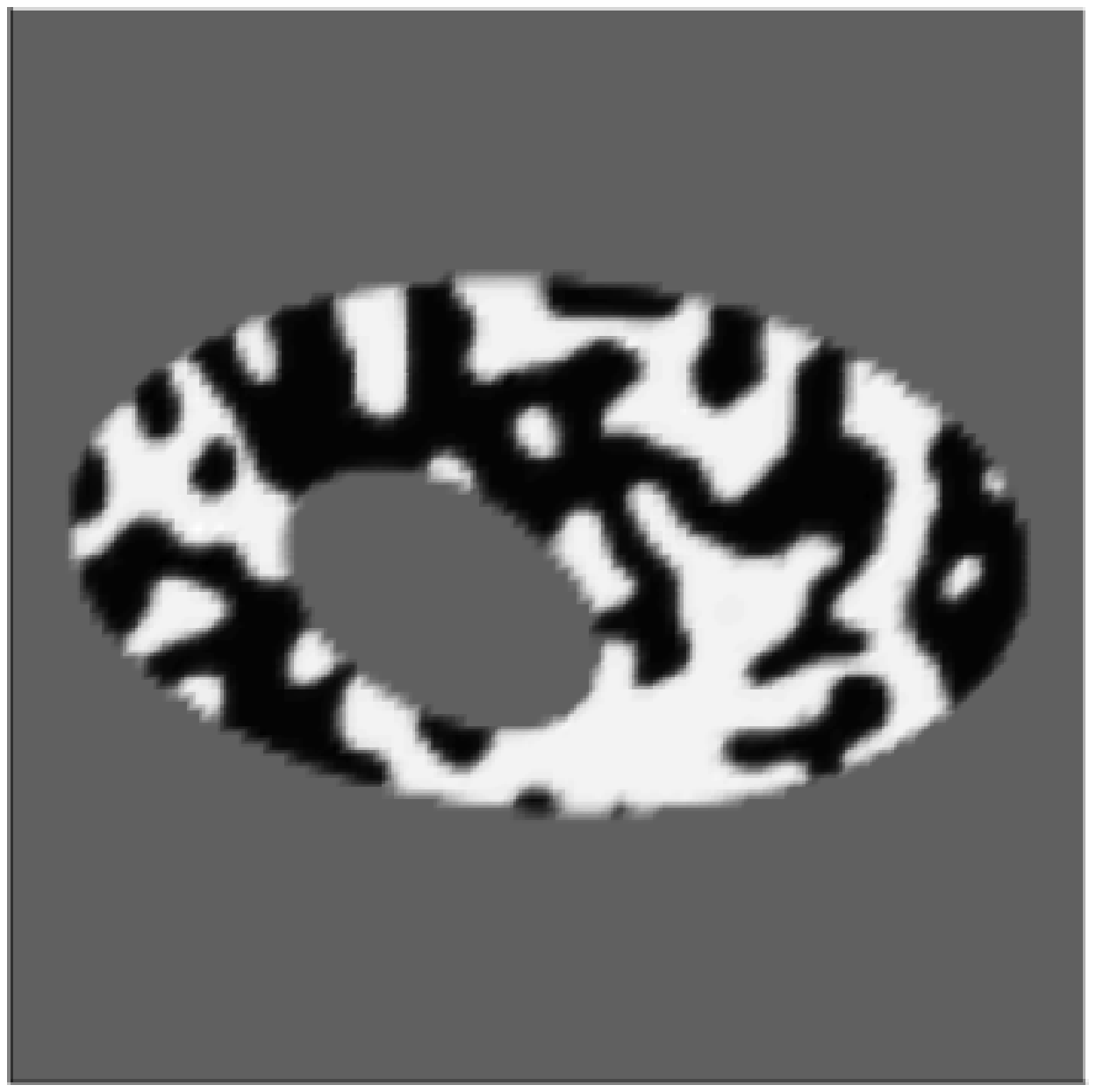} &
\includegraphics[height=1.5in]{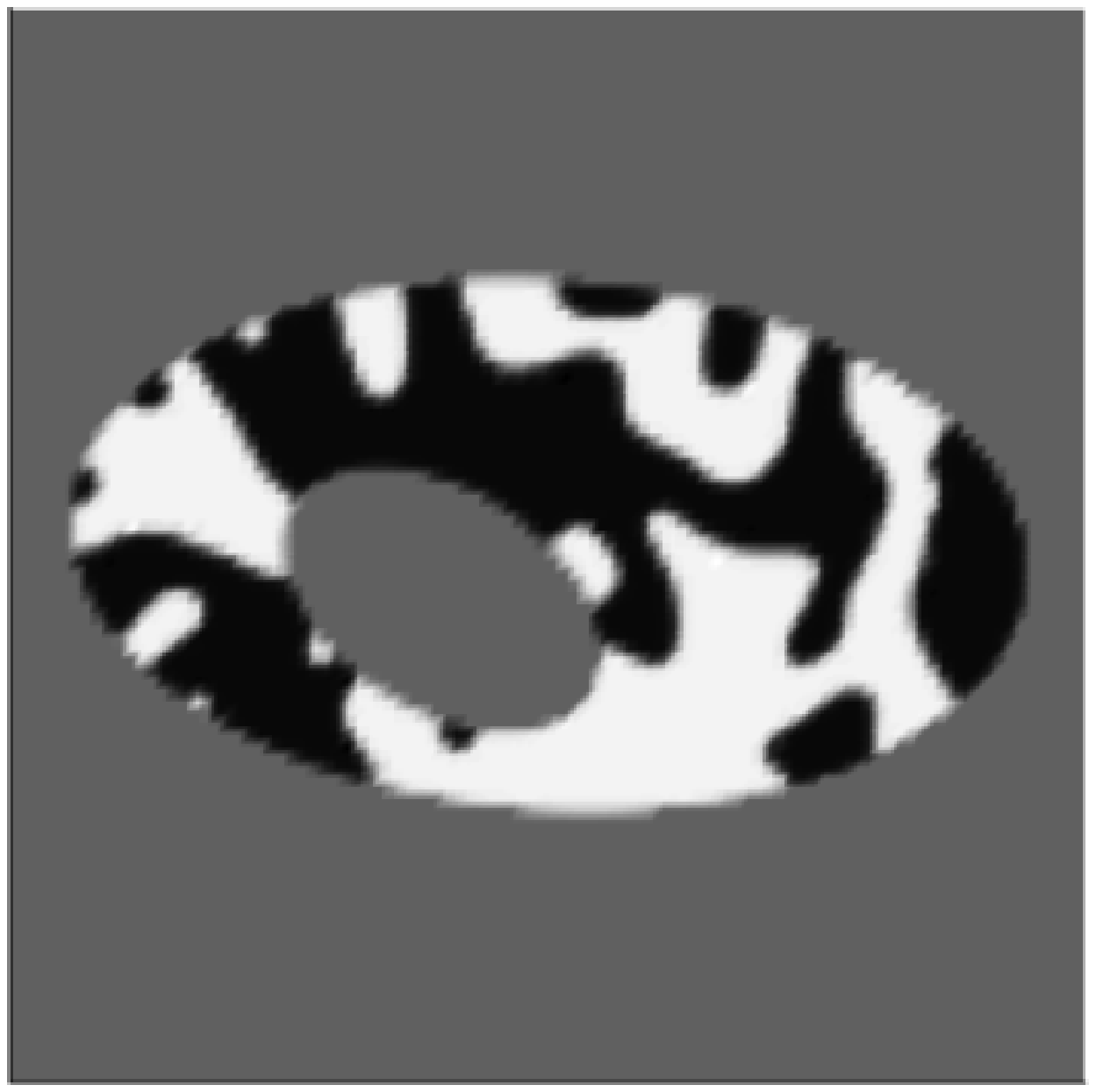} 
\end{array}$
  \caption{\em The solution to Example 3, 100 time steps of size one are taken.  The first picture shows the initial conditions, and the final picture is at $t=100$. }
  \label{figure4}
\end{figure}

{\it{EXAMPLE 3.}} 
In this example, we demonstrate that our methods can be applied to much more complex equations. 
Here, we solve the Allen-Cahn equation 
\begin{align*}
   u_t - \epsilon \Delta u & =  u(1-u^2), &   \x \in \Omega, \; t>0, \\
   u(\x,t) & =   0, & \x \in \Gamma_0  \\
   u(\x,t) & =  0,  &  \x \in \Gamma_1,
\end{align*}
where $\epsilon = 10^{-5}$. 
We initialize the solution with random values uniformly distributed on $[-\frac{1}{2}, \frac{1}{2}]$.  
The general behaviour of solutions to the Allen-Cahn equation is well known: 
the stable stationary solutions are $u=1$ and $u=-1$ and the solution exhibits coarsening towards these values.
The presence of physical boundaries can create more complex patterns, as seen in \figr{figure4}. 

\section{Conclusions}
We have presented an investigation on coupling available fast algorithms for integral equations  in order to investigate  problems of scientific interest.  
We have shown that for time-dependent problems, doing a temporal discretization first, and then reformulating the resulting elliptic equation as an integral equation, appears to have significant potential for handling complicated time dependent problems on complex domains. 
This paper is meant to be a preliminary investigation only; more work still has to be done to make this approach amenable to large-scale problems with a variety of boundary conditions. 
Specifically, we need to develop a more robust and general purpose method of extending functions defined on a complex domain onto a regular domain in a way that does not result in excess grid refinement in the vicinity of physical boundaries. 
Future applications will likely include developing integral equation methods for the incompressible Navier-Stokes and Euler equations.

\noindent
{\bf Acknowledgements:}
We wish to thank the organizers of the conference in {\it Advances in Boundary Integral Equations and Related Topics} for their invitation and the opportunity to participate in this special issue. We also wish to thank Rustum Choksi and Yves Van Gennip for their help with Example 3, and Jingfang Huang for his considerable help with his volume potential solver. 

\bibliographystyle{model1b-num-names}

\end{document}